\documentclass[12pt]{amsart}

\usepackage[cp1251]{inputenc}
\usepackage[T2A]{fontenc}
\usepackage{amsfonts}
\usepackage{amsbsy}
\usepackage{amssymb}
\usepackage[english]{babel}
\usepackage{hyperref}
\usepackage{wrapfig}

\catcode`~=11 
\newcommand{\urltilde}{\kern -.15em\lower .7ex\hbox{~}\kern .04em}
\catcode`~=13 

\renewcommand{\abovecaptionskip}{0pt}
\renewcommand{\belowcaptionskip}{6pt}
\makeatletter

\renewcommand{\@makecaption}[2]{
\vspace{\abovecaptionskip}%
\sbox{\@tempboxa}{#1 #2}%
\global\@minipagefalse \hbox to \hsize {{\scshape \hfil #1 #2\hfil}}
\vspace{\belowcaptionskip}}

\makeatother

\newcommand{\mf}{\mathfrak}

\newcommand{\mr}{\mathrm}

\newcommand{\reg}{\mr{reg}}

\renewcommand{\ge}{\geqslant}

\newcommand{\Ker}{\operatorname{Ker}}

\newcommand{\SL}{\operatorname{SL}}
\newcommand{\GL}{\operatorname{GL}}

\newcommand{\diag}{\operatorname{diag}}

\DeclareMathOperator{\otimesZ}{\otimes\hspace{1pt}\rule[-3pt]{0pt}{0pt}_{\mathbb{Z}}}
\DeclareMathOperator{\Ad}{\mr{Ad}}
\DeclareMathOperator{\Supp}{\mr{Supp}}
\DeclareMathOperator{\hgt}{ht} 

\newtheorem{theorem}{Theorem}

\newtheorem{proposition}{Proposition}
\newtheorem{lemma}{Lemma}
\newtheorem{corollary}{Corollary}

\newtheorem*{question*}{Question}

\theoremstyle{definition}

\newtheorem{definition}{Definition}
\newtheorem*{definition*}{Definition}
\newtheorem{example}{Example}

\theoremstyle{remark}

\newtheorem{remark}{Remark}

\begin{document}

\renewcommand{\proofname}{Proof}
\renewcommand{\abstractname}{Abstract}
\renewcommand{\refname}{References}
\renewcommand{\figurename}{Figure}
\renewcommand{\tablename}{Table}

\title[Normalizers of solvable spherical subgroups]{Normalizers of solvable spherical subgroups}

\author{Roman Avdeev}

\thanks{Partially supported by Russian Foundation for Basic Research (grant no. 09-01-00648) and by Alexander von Humboldt foundation}

\address{Chair of Higher Algebra, Department of Mechanics and Mathematics,
Moscow State University, 1, Leninskie Gory, Moscow, 119992, Russia}

\email{suselr@yandex.ru}


\subjclass[2010]{20G07, 14M27}

\keywords{Algebraic group, spherical subgroup, solvable subgroup,
normalizer}

\begin{abstract}
For an arbitrary connected solvable spherical subgroup~$H$ of a
connected semisimple algebraic group~$G$ we compute the group
$N_G(H)$, the normalizer of~$H$ in~$G$. Thereby we complete a
classification of all (not necessarily connected) solvable spherical
subgroups in semisimple algebraic groups.
\end{abstract}

\maketitle

\sloppy

\section{Introduction} \label{sect_intro}

Let $G$ be a connected semisimple complex algebraic group. A closed
subgroup $H \subset G$ (resp. a homogeneous space~$G/H$) is said to
be \textit{spherical} if a Borel subgroup $B \subset G$ has an open
orbit in~$G/H$. The latter condition is equivalent to the relation
$\mathfrak b + (\Ad g)\mathfrak h = \mathfrak g$ holding for some
element $g \in G$, where $\mathfrak b, \mathfrak h, \mathfrak g$ are
the Lie algebras of the groups $B, H, G$, respectively, and $\Ad
\colon G \to \GL(\mathfrak g)$ is the adjoint representation of~$G$.
In particular, this implies that the sphericity of $H$ (as~well as
of~$G/H$) is a \textit{local} property, that is, it depends only on
the pair of algebras~$(\mathfrak g,\mathfrak h)$. Thus one may talk
about spherical subalgebras in~$\mathfrak g$. Equivalently, a
subgroup $H$ is spherical if and only if its connected component of
the identity $H^0$ is spherical.

As the sphericity is a local property, the classification of
spherical subgroups in~$G$ can be naturally divided into two stages.
The first stage is the classification of \textit{connected}
spherical subgroups in~$G$, which is equivalent to the
classification of spherical subalgebras in~$\mathfrak g$. At the
second stage, for every connected spherical subgroup $H_0 \subset G$
one has to find all subgroups $H \subset G$ with $H^0 = H_0$.

By the present time, within the first stage of the classification an
explicit description (in terms of Lie algebras) of all spherical
subgroups has been obtained inside two `opposite' classes of
subgroups, namely, reductive and solvable. For simple groups, the
list of connected reductive spherical subgroups was obtained
in~\cite{Kr}, and that for non-simple semisimple groups was obtained
in~\cite{Mi} and, independently, in~\cite{Br} (see also~\cite{Yak}
for a more accurate formulation). Connected solvable spherical
subgroups were described in~\cite{Avd_solv}.

Let $H_0$ be a connected subgroup of~$G$. Then every subgroup $H
\subset G$ with $H^0 = H_0$ is contained in the subgroup $N_G(H_0)$,
the normalizer of $H_0$ in~$G$. Thereby one establishes a one-to-one
correspondence between finite extensions of the subgroup $H_0$ and
finite subgroups of the group~$N_G(H_0)/H_0$. Thus the second stage
of the classification reduces to the computation of the
group~$N_G(H_0)$ for every connected spherical subgroup~$H_0$.

In~\cite{Avd_exc}, the normalizers of most of reductive spherical
subgroups of simple groups were found. In the present paper, we find
the normalizers of all connected solvable spherical subgroups of
semisimple algebraic groups. For this purpose, we use the structure
theory of connected solvable spherical subgroups developed
in~\cite{Avd_solv}. The main results of this paper are
Theorem~\ref{thm_main_result} and Proposition~\ref{prop_root_in_P}
(see \S\,\ref{sect_result}).

It is known that for an arbitrary spherical subgroup $H \subset G$
the group $N_G(H)/H$ is commutative (and even diagonalizable,
see~\cite[Corollary~5.2]{BrP}). Therefore:

(1) a spherical subgroup $H \subset G$ is solvable if and only if
the subgroup~$H^0$ is solvable; this means that the description of
all finite extensions of connected solvable spherical subgroups
provides a classification of all (not necessarily connected)
solvable spherical subgroups in semisimple algebraic groups;

(2) the subgroup $N_G(H)$ coincides with the subgroup $N_G(H^0)$; in
view of this the results of this paper enable one to compute the
normalizers of arbitrary solvable spherical subgroups in semisimple
algebraic groups.

\begin{remark}
By this moment, a program of classification of \textit{all}
spherical homogeneous spaces in combinatorial terms (so-called
\textit{spherical homogeneous data}) is completed. This program was
initiated by Luna in~\cite{Lu01} (more details see
in~\cite[\S\,30.11]{Tim}). However this classification is implicit
in the following sense: so far there is no general procedure of
constructing a spherical subgroup given by its set of combinatorial
data.
\end{remark}

This paper is organized as follows. In~\S\,\ref{sect_notation} we
list some notation and conventions used in the paper. In
\S\,\ref{sect_preliminaries} we collect initial facts from the
structure theory of connected solvable spherical subgroups that are
needed for formulation of the main results. In \S\,\ref{sect_result}
we state the main results of the paper
(Theorem~\ref{thm_main_result} and Proposition~\ref{prop_root_in_P})
and give some examples. In \S\,\ref{sect_further_results} we collect
further facts from the structure theory of connected solvable
spherical subgroups that are needed for proving the main results. At
last, the proofs of the main results are given in
\S\,\ref{sect_proof}.

This work was completed during the author's stay at hospitable
Bielefeld University in June and July of 2011. This stay was
supported by Alexander von Humboldt Foundation.

The author expresses his gratitude to E.\,B.~Vinberg for reading the
paper and valuable comments.

\section{Some notation and conventions}
\label{sect_notation}

In this paper the base field is the field $\mathbb C$ of complex
numbers. All groups are assumed to be algebraic and their subgroups
closed in the Zarisky topology. The Lie algebras of groups denoted
by capital Latin letters are denoted by the corresponding small
German letters.

Some notation:

$e$ is the identity element of an arbitrary group;

$|X|$ is the cardinality of a finite set~$X$;

$L^0$ is the connected component of the identity of a group~$L$;

$\mf X(L)$ is the character lattice of a group~$L$;

$N_L(K)$ is the normalizer of a subgroup $K$ in a group~$L$;

$G$ is an arbitrary connected semisimple algebraic group;

$B \subset G$ is a fixed Borel subgroup of~$G$;

$T \subset B$ is a fixed maximal torus of~$G$;

$U \subset B$ is the maximal unipotent subgroup of~$G$ contained
in~$B$;

$W = N_G(T)/T$ is the Weyl group of $G$ with respect to~$T$;

$\theta \colon N_G(T) \to W$ is the canonical homomorphism;

$Q = \mathfrak X(T)\otimesZ \mathbb Q$ is the rational vector space
spanned by~$\mathfrak X(T)$;

$\langle A \rangle$ is the linear span in $Q$ of a subset $A \subset
\mathfrak X(T)$;

$(\cdot\,, \cdot)$ is a fixed inner product in $Q$ invariant
under~$W$;

$\Delta \subset \mathfrak X(T)$ is the root system of $G$ with
respect to~$T$;

$\Delta_+ \subset \Delta$ is the subset of positive roots with
respect to~$B$;

$\Pi \subset \Delta_+$ is the set of simple roots;

$r_\alpha \in W$ is the simple reflection corresponding to a root
$\alpha \in \Pi$;

$\mathfrak g_\alpha \subset \mathfrak g$ is the root subspace
corresponding to a root $\alpha \in \Delta$;

$\Sigma(\widetilde \Pi)$ is the Dynkin diagram of a subset
$\widetilde \Pi \subset \Pi$.

\vspace{1em}

For every root $\alpha = \sum\limits_{\delta \in \Pi}k_\delta \delta
\in \Delta_+$ we consider its \emph{support} $\Supp \alpha =
\{\delta \mid k_\delta > 0\}$ and \emph{height} $\hgt \alpha =
\sum\limits_{\delta \in \Pi}k_\delta$. If $\alpha \in \Delta_+$, we
put $\Delta(\alpha) = \Delta \cap \langle \Supp \alpha \rangle$ and
$\Delta_+(\alpha) = \Delta_+ \cap \langle \Supp\alpha \rangle$. The
set $\Delta(\alpha)$ is an indecomposable root system whose set of
simple roots is $\Supp \alpha$; the set of positive roots of this
root system coincides with $\Delta_+(\alpha)$.

Let $L$ be a group and let $L_1, L_2$ be subgroups of it. We write
$L = L_1 \rightthreetimes L_2$ if $L$ is a semidirect product of
$L_1, L_2$, that is, $L = L_1L_2$, $L_1 \cap L_2 = \{e\}$, and $L_2$
is a normal subgroup of~$L$.

By abuse of language, we identify roots in $\Pi$ and the
corresponding nodes of the Dynkin diagram of~$\Pi$.

When we say that two nodes of a Dynkin diagram are joined by and
edge we mean that the edge may be multiple.

The enumeration of nodes (that is, of simple roots) of connected
Dynkin diagrams is the same as in the book~\cite{VO}.

\section{Preliminaries on connected solvable spherical subgroups} \label{sect_preliminaries}

In this subsection we collect basic definitions and initial facts of
the structure theory of connected solvable spherical subgroups in
semisimple algebraic groups (see~\cite{Avd_solv}) needed for
formulation of the main results of this paper.

Let $H \subset B$ be a connected solvable subgroup and $N \subset U$
its unipotent radical. We say that the subgroup $H$ is
\textit{standardly embedded in~$B$} (with respect to~$T$) if the
subgroup $S = H \cap T \subset T$ is a maximal torus of~$H$.
Evidently, in this situation we have $H = S \rightthreetimes N$.
Every connected solvable subgroup of $G$ is conjugate to a subgroup
that is standardly embedded in~$B$.

Suppose that a connected solvable subgroup ${H \subset G}$
standardly embedded in $B$ is fixed. As above, we put $S = H \cap T$
and $N = H \cap U$ so that $H = S \rightthreetimes N$. We denote by
$\tau \colon \mathfrak X(T)\to\mathfrak X(S)$ the character
restriction map from $T$ to~$S$. Let $\Phi = \tau(\Delta_+) \subset
\mathfrak X(S)$ be the weight system of the action of $S$ on
$\mathfrak u$ by means of the adjoint representation of~$G$. We have
$\mathfrak u = \bigoplus \limits_{\varphi \in \Phi} \mathfrak
u_\varphi$, where $\mathfrak u_\varphi \subset \mathfrak u$ is the
weight subspace of weight $\varphi$ with respect to~$S$. For
$\varphi \in \Phi$ we put $\mathfrak n_\varphi = \mathfrak n \cap
\mathfrak u_\varphi$. Evidently, we have $\mathfrak n = \bigoplus
\limits_{\varphi \in \Phi} \mathfrak n_\varphi$. For every $\varphi
\in \Phi$ let $c_\varphi$ denote the codimension of the subspace
$\mathfrak n_\varphi$ in the space $\mathfrak u_\varphi$.

In the notation introduced above, the following criterion of
sphericity of $H$ takes place.

\begin{theorem}[{\rm \cite[Theorem~1]{Avd_solv}}]\label{thm_solv_spher}
The following conditions are equivalent:

\textup{(1)} $H$ is spherical in~$G$;

\textup{(2)} $c_\varphi \leqslant 1$ for every $\varphi \in \Phi$,
and all weights with $c_\varphi = 1$ are linearly independent
in~$\mathfrak X(S)$.
\end{theorem}

Later on, we assume that $H \subset G$ is a connected solvable
spherical subgroup standardly embedded in~$B$ and preserve all the
notations introduced above. Put $\Psi = {\{\alpha \in \Delta_+ \mid
\mathfrak g_\alpha \not \subset \mathfrak n\}} \subset \Delta_+$.

\begin{definition}
The roots in the set $\Psi$ are said to be \textit{active}.
\end{definition}

Let $\varphi_1, \ldots, \varphi_m$ denote the weights $\varphi \in
\Phi$ with $c_\varphi = 1$. For $i = 1, \ldots, m$ we put $\Psi_i =
\{\alpha \in \Psi \mid \tau(\alpha) = \varphi_i\}$. It is clear that
$\Psi = \Psi_1 \cup \Psi_2 \cup \ldots \cup \Psi_m$ and ${\Psi_i
\cap \Psi_j = \varnothing}$ for $i \ne j$. A~key role in the
structure theory of connected solvable spherical subgroups is played
by the following proposition.

\begin{proposition}[{\rm\cite[Proposition~1]{Avd_solv}}] \label{prop_Psi's}
Suppose that $1 \leqslant i,j \leqslant m$ and different roots
$\alpha \in \Psi_i$, $\beta \in \Psi_j$ are such that $\gamma =
\beta - \alpha \in \Delta_+$. Then $\Psi_i + \gamma \subset \Psi_j$.
\end{proposition}

In particular, Proposition~\ref{prop_Psi's} is used in the proof of
the following theorem.

\begin{theorem}[{\rm \cite[Theorem~2]{Avd_solv}}]
\label{thm_approximation} Up to conjugation by an element of~$T$,
the subgroup $H$ is uniquely determined by the pair $(S, \Psi)$.
\end{theorem}

The set of active roots has the following property
(see~\cite[Lemma~4]{Avd_solv}): if $\alpha$ is an active root and
$\alpha = \beta + \gamma$ for some roots $\beta, \gamma \in
\Delta_+$, then exactly one of the two roots $\beta, \gamma$ is
active. Taking this property into account, we say that an active
root $\beta$ is \textit{subordinate} to an active root~$\alpha$ if
$\alpha = \beta + \gamma$ for some root $\gamma \in \Delta_+$. For
every active root $\alpha$ we denote by $F(\alpha)$ the set
consisting of $\alpha$ and all active roots subordinate to~$\alpha$.
An active root $\alpha$ is said to be \textit{maximal} if it is not
subordinate to any other active root. Let $\mathrm M$ denote the set
of maximal active roots.

Let us mention the following corollary from
Proposition~\ref{prop_Psi's}.

\begin{corollary} \label{crl_either_or}
For every $i = 1, \ldots, m$, we have either $\Psi_i \subset \mathrm
M$ or $\Psi_i \cap \mathrm M = \varnothing$.
\end{corollary}

\begin{proposition}[{\rm \cite[Proposition~3]{Avd_solv}}]
\label{prop_assoc_root} Let $\alpha$ be an active root. Then there
exists a unique simple root \mbox{$\pi(\alpha) \in \Supp \alpha$}
with the following property: if $\alpha = \beta + \gamma$ for some
roots $\beta, \gamma \in \Delta_+$, then the root $\beta$ is active
if and only if $\pi(\alpha) \notin \Supp \beta$ \textup{(}and so the
root $\gamma$ is active if and only if $\pi(\alpha) \notin \Supp
\gamma$\textup{)}.
\end{proposition}

From Proposition~\ref{prop_assoc_root} it follows that for every
active root $\alpha$ the set $F(\alpha)$ is uniquely determined by
the simple root~$\pi(\alpha)$. Since every active root is
subordinate to a maximal active root
(see~Proposition~\ref{prop_F&Supp}(d)
in~\S\,\ref{sect_further_results}
or~\cite[Corollary~3(b)]{Avd_solv}), we obtain that the whole set
$\Psi$ is uniquely determined by the subset $\mathrm M$ and the
map~$\pi \colon \mathrm M \to \Pi$.

We now introduce an equivalence relation on $\mathrm M$ as follows.
For any two roots $\alpha, \beta \in \mathrm M$ we write $\alpha
\sim \beta$ if and only if $\tau(\alpha) = \tau(\beta)$.

To the subgroup $H$ we assign the set $\Upsilon_0(H) = (\mathrm M,
\pi, \sim)$, where $\pi$ is regarded as a map from $\mathrm M$
to~$\Pi$.

\begin{proposition}[{\rm \cite[Remark~4]{Avd_solv}}]
\label{prop_unip_rad} Up to conjugation by an element in~$T$, the
unipotent radical $N$ of $H$ is uniquely determined by the
set~$\Upsilon_0(H)$.
\end{proposition}

\section{Formulation of the main results} \label{sect_result}

In this subsection we state the main results of this paper
(Theorem~\ref{thm_main_result} and Proposition~\ref{prop_root_in_P})
and give some examples.

Let $H \subset G$ be a connected solvable spherical subgroup
standardly embedded in~$B$. We preserve all the notations introduced
in~\S\,\ref{sect_preliminaries}.

For a fixed $i \in \lbrace 1, \ldots, m \rbrace$, the subspace
$\mathfrak n_{\varphi_i}$ of the space $\mathfrak u_{\varphi_i}$ is
defined by the vanishing of a linear function~$\xi_i$, which is
determined up to proportionality. It follows from the definition of
the set $\Psi_i$ that for each $\alpha \in \Psi_i$ the restriction
of $\xi_i$ to the subspace $\mathfrak g_{\alpha}$ does not vanish.

Let $L$ denote the sublattice in~$\mathfrak X(T)$ generated by all
elements of the form $\alpha - \beta$, where $\alpha, \beta \in
\mathrm M$ and $\tau(\alpha) = \tau(\beta)$. In view of
Proposition~\ref{prop_Psi's} and Corollary~\ref{crl_either_or}, $L$
is also generated by all elements of the form $\alpha - \beta$,
where $\alpha, \beta \in \mathrm \Psi$ and $\tau(\alpha) =
\tau(\beta)$. Further, we introduce the lattice $L_0 = \langle L
\rangle \cap \mathfrak X(T) \supset L$ and denote by~$A$ (resp.
by~$A_0$) the subgrooup in~$T$ defined by the condition that all
characters in~$L$ (resp. in~$L_0$) are equal to one. Since the
sublattice $L_0 \subset \mathfrak X(T)$ is primitive, we have ${A^0
= A_0}$ and $A/A_0 \simeq L_0/L$. From what we said in the previous
paragraph it is clear that the group~$A$ (resp.~$A_0$) is the
largest subgroup (resp. the largest connected subgroup) in~$T$
normalizing~$N$.

\begin{definition} \label{dfn_regular}
An active root $\delta$ is said to be \textit{regular} if the
projection of the subspace $\mathfrak n \subset \mathfrak u$ to the
root subspace $\mathfrak g_\delta$ along the sum of the other root
subspaces is zero.
\end{definition}

Let $\delta$ be an active root and choose $i \in \{1, \ldots, m\}$
such that $\delta \in \Psi_i$. Then from
Definition~\ref{dfn_regular} it follows that the root $\delta$ is
regular if and only if $|\Psi_i| = 1$, that is, $\Psi_i =
\{\delta\}$.

Let $\Psi^{\reg} \subset \Psi$ be the set of regular active roots.

We denote by $\mathrm P$ the set of active roots~$\alpha$ satisfying
the following two conditions:

(1) $\alpha \in \Psi^{\reg} \cap \Pi$;

(2) $(\alpha, \beta) = 0$ for every root $\beta \in \Psi \backslash
\{\alpha\}$.

It is easy to see that for every root $\delta \in \mathrm P$ and
every element $\rho \in \theta^{-1}(r_\delta)$ the subgroup $H' =
\rho H \rho^{-1}$ is also standardly embedded in~$B$, at that,
$\Upsilon_0(H') = \Upsilon_0(H)$ in view of condition~(2). By
Proposition~\ref{prop_unip_rad} we have $N' = tNt^{-1}$ for some
element $t \in T$. This implies that $t^{-1}\rho \in N_G(N)$, whence
$\theta^{-1}(r_\delta) \cap N_G(N) \ne \varnothing$.

For every root $\delta \in \mathrm P$ we fix an arbitrary element
$\rho_\delta \in \theta^{-1}(r_\delta) \cap N_G(N)$. We put $\mathrm
P_S = \{\delta \in \mathrm P \mid r_\delta(\Ker \tau) = \Ker \tau\}
= \{\delta \in \mathrm P \mid \rho_\delta S \rho_\delta^{-1} = S\}$.

We now are able to state the main results of this paper.

\begin{theorem} \label{thm_main_result}
The group $N_G(H)$ is generated by the groups $A$, $N$ and all
elements~$\rho_\delta$, where $\delta$ runs over the set $\mathrm
P_S$. In particular, $N_G(H)^0 = A_0 \rightthreetimes N$.
\end{theorem}

\begin{corollary}
There are isomorphisms $N_G(H)/H \simeq A/S \times (\mathbb Z/2
\mathbb Z)^r$ and $N_G(H)/N_G(H)^0  \simeq A/A_0 \times (\mathbb Z/2
\mathbb Z)^r \simeq L_0/L \times (\mathbb Z/2 \mathbb Z)^r$, where
$r = |\mathrm P_S|$.
\end{corollary}

\begin{remark}
The group $N_G(H)^0$ was computed earlier
in~\cite[\S\,5.5]{Avd_solv}.
\end{remark}

Let $\alpha$ be an active root. A simple root $\delta \in \Supp
\alpha$ is said to be \textit{terminal with respect to $\Supp
\alpha$} if in the diagram $\Sigma(\Supp \alpha)$ the node
corresponding to~$\delta$ is joined by an edge with exactly one
node.

The following proposition enables one to find the set~$\mathrm P$
explicitly given the set $\Upsilon_0(H) = (\mathrm M, \pi, \sim)$.

\begin{proposition} \label{prop_root_in_P}
A root $\alpha \in \Pi$ is contained in the set $\mathrm P$ in
exactly one of the following two cases:

Case~\textup{1}. The following conditions are fulfilled:

{\addtolength{\leftskip}{\parindent}

\textup{(1)} $\alpha \in \mathrm M$;

\textup{(2)} $\alpha \not\sim \beta$ for every root $\beta \in
\mathrm M \backslash \{\alpha\}$;

\textup{(3)} for every root $\beta \in \mathrm M\backslash
\{\alpha\}$ the diagram $\Sigma(\{\alpha\} \cup \Supp \beta)$ is
disconnected.

}

Case~\textup{2}. $\alpha \notin \mathrm M$ and there is a root
$\beta \in \mathrm M \backslash \{\alpha\}$ with the following
properties:

{\addtolength{\leftskip}{\parindent}

\textup{(1)} $\beta = \sum \limits_{\delta \in \Supp \beta} \delta$;

\textup{(2)} $\alpha \in \Supp \beta$ and $\alpha$ is terminal with
respect to $\Supp \beta$;

\textup{(3)} $\pi(\beta) \ne \alpha$;

\textup{(4)} there is a root $\alpha' \in \Supp \beta$ such that in
the diagram $\Sigma(\Supp \beta)$ the nodes $\alpha$ and $\alpha'$
are joined by a double edge with the arrow directed to~$\alpha$;

\textup{(5)} for every root $\gamma \in \mathrm M \backslash
\{\beta\}$ the diagram $\Sigma(\{\alpha\} \cup \Supp \gamma)$ is
disconnected.

}
\end{proposition}

\begin{example} \label{ex_adjoint}
Let $G$ be adjoint (that is, its center is trivial). It is well
known that in this case the lattice $\mathfrak X(T)$ is generated by
the set~$\Pi$. Let us show that $A = A_0$. Replacing $H$ by a
subgroup conjugate by an appropriate element in $N_G(T)$, we may
assume that each root in $\mathrm M$ equals the sum of all simple
roots in its support (see~\cite[Proposition~17]{Avd_solv}). Then
condition~$(\mathrm C)$ (see Theorem~\ref{thm_classification} in
\S\,\ref{sect_further_results}) implies that the set $\mathrm M$
generates a primitive sublattice in $\mathfrak X(T)$. It follows
that the sublattice $L$ is also primitive, whence $L = L_0$ and $A =
A_0$. Hence $N_G(H)/N_G(H)^0 \simeq (\mathbb Z/2\mathbb Z)^r$, where
$r = |\mathrm P_S|$.
\end{example}

\begin{remark}
The equality $A = A_0$ in case of an adjoint group~$G$ implies that
$A = Z(G)A_0$ in case of an arbitrary group~$G$.
\end{remark}

\begin{example} \label{ex_sober}
Suppose that $S = A_0$. Then $H = N_G(H)^0$ and $\Ker \tau = L_0$.
It is easy to see that for every root $\delta \in \mathrm P$ the
reflection~$r_\delta$ preserves the subspace $\langle L_0 \rangle
\subset Q$ and acts trivially on it. This implies that $\mathrm P_S
= \mathrm P$. Therefore $N_G(H)/N_G(H)^0 \simeq L_0/L \times
(\mathbb Z/ 2\mathbb Z)^r$, where $r = |\mathrm P|$.
\end{example}

\begin{example}
Suppose that $G$ is simple and adjoint and $H$ satisfies the
condition $\bigcup \limits_{\alpha \in \mathrm M} \Supp \alpha =
\Pi$. Then by condition~$(\mathrm C)$ (see
Theorem~\ref{thm_classification} in \S\,\ref{sect_further_results})
and Example~\ref{ex_adjoint} we obtain $S = A_0$. In view of
Proposition~\ref{prop_root_in_P} there are the following
possibilities depending on the type of the root system~$\Delta$:

(1) if $\Delta$ is of type different from $\mathsf A_1$ or $\mathsf
B_n$ ($n \ge 2$), then $\mathrm P = \varnothing$ and $N_G(H) = H$;

(2) if $\Delta$ is of type $\mathsf A_1$, then, evidently, $H = T$,
$\mathrm P = \Pi$, and $N_G(H)/H \simeq \mathbb Z/ 2\mathbb Z$;

(3) if $\Delta$ is of type $\mathsf B_n$ ($n \ge 2$), then either
$\mathrm P = \varnothing$ (which implies $N_G(H) = H$) or $|\mathrm
P| = 1$ (which implies $N_G(H)/H \simeq \mathbb Z/ 2\mathbb Z$).
Suppose that $\Pi = \{\alpha_1, \ldots, \alpha_n\}$. Then as an
example of the first situation we can take $\mathrm M = \{\beta\}$,
where $\beta = \alpha_1 + \ldots + \alpha_n$, and $\pi(\beta) =
\alpha_n$. As an example of the second situation we can take
$\mathrm M = \{\beta\}$, where $\beta = \alpha_1 + \ldots +
\alpha_n$, and $\pi(\beta) \ne \alpha_n$; we have $\mr P = \lbrace
\alpha_n \rbrace$.
\end{example}

\begin{example} \label{ex_counterex}
Suppose that $G = \SL_3$ and the groups $B, U, T$ consist of all
upper-triangular, upper unitriangular, diagonal matrices,
respectively, contained in~$G$. For $t = \diag(t_1,t_2,t_3) \in T$
and $k = 1, 2$ we put $\alpha_k(t) = t_kt_{k+1}^{-1}$. Let $H$ be
the group of matrices of the form
\begin{equation} \label{eqn_matrix}
\begin{pmatrix}
t & 0 & a\\
0 & 1 & b\\
0 & 0 & t^{-1}
\end{pmatrix},
\end{equation}
where $a, b, t \ne 0$ are arbitrary numbers. Then $S = \{\diag(t, 1,
t^{-1}) \mid t \ne 0\}$ and the subgroup $N$ consists of all
matrices of the form~(\ref{eqn_matrix}) with $t = 1$. We have $A =
A_0 = T$ and $\Psi = \mathrm M = \mathrm P = \{\alpha_1\}$. Besides,
$\langle \Ker \tau \rangle = \langle \alpha_1 - \alpha_2 \rangle$,
whence $r_{\alpha_1} \langle \Ker \tau \rangle \ne \langle \Ker \tau
\rangle$ and $\mathrm P_S = \varnothing$. Thus, $N_G(H) = T
\rightthreetimes N$. Further, we have $N_G(N_G(H)) = N_G(H) \cup
\rho N_G(H)$, where
$$
\rho =
\begin{pmatrix}
0 & 1 & 0\\
1 & 0 & 0\\
0 & 0 & -1
\end{pmatrix}.
$$
In particular, we find that $N_G(N_G(H)) \ne N_G(H)$.
\end{example}

\begin{remark}
Example~\ref{ex_counterex} disproves Theorem~4.3(iii) in~\cite{Br97}
and Lemma~30.2 in~\cite{Tim}, which claim that $N_G(N_G(H)) =
N_G(H)$ for an arbitrary spherical subgroup $H \subset G$.
Nevertheless, this error does not influence the truth of other
results of the general theory of spherical homogeneous spaces.
\end{remark}

\begin{remark}
Applying Theorem~\ref{thm_main_result} to the group $N_G(H)^0 = A_0
\rightthreetimes N$ (see Example~\ref{ex_sober}), we obtain that the
group $N_G(N_G(H))$ is generated by the groups $A$, $N$, and all
elements $\rho_\delta$, where $\delta \in \mr P$. In particular, the
equality $N_G(N_G(H)) = N_G(H)$ holds if and only if $\mathrm P =
\mathrm P_S$.
\end{remark}

\section{Further results on connected solvable spherical subgroups} \label{sect_further_results}

In this subsection we collect all results from the structure theory
and classification of connected solvable spherical subgroups needed
for proving Theorem~\ref{thm_main_result} and
Proposition~\ref{prop_root_in_P}.

Throughout this subsection we denote by $H$ and $H'$ two (not
necessarily different) connected solvable spherical subgroups
standardly embedded in~$B$. The notations $S$, $N$, $\Psi$, $\pi$,
\ldots (resp. $S'$, $N'$, $\Psi'$, $\pi'$, \ldots) refer to the
subgroup~$H$ (resp. to the subgroup~$H'$).

\begin{proposition}[{\rm \cite[Lemma~7]{Avd_solv}}] \label{prop_F&Supp}
Suppose that $\alpha \in \Psi$. Then:

\textup{(a)} $|F(\alpha)| = |\Supp \alpha|$;

\textup{(b)} all weights $\tau(\beta)$, where $\beta \in \Supp
\alpha$, are linearly independent in~$\mathfrak X(S)$;

\textup{(c)} $\langle F(\alpha) \rangle = \langle \Supp \alpha
\rangle$;

\textup{(d)} if $\beta \in \Psi$ and $\Supp \beta \subset \Supp
\alpha$, then $\beta \in F(\alpha)$.
\end{proposition}

\begin{proposition}[{\rm\cite[Corollary~6]{Avd_solv}}]
\label{prop_bijection} Suppose that $\alpha \in \Psi$. Then the map
$\pi \colon {F(\alpha) \to \Supp \alpha}$ is a bijection.
\end{proposition}

\begin{proposition}[{\rm\cite[Lemma~10]{Avd_solv}}] \label{prop_equal_pi}
Suppose that roots $\alpha, \beta \in \Psi$ are such that
$\pi(\alpha) = \pi(\beta)$. Then $\tau(\alpha) = \tau(\beta)$.
\end{proposition}

\begin{proposition} \label{prop_not_contained}
Suppose that different roots $\alpha, \beta \in \Psi$ are such that
$\tau(\alpha) = \tau(\beta)$ and $\pi(\alpha) = \pi(\beta)$. Then
there are roots $\widetilde \alpha \in F(\alpha)$ and $\widetilde
\beta \in F(\beta)$ such that $\tau(\widetilde \alpha) =
\tau(\widetilde \beta)$, ${\pi(\widetilde \alpha) \notin \Supp
\beta}$, and $\pi(\widetilde \beta) \notin \Supp \alpha$.
\end{proposition}

This proposition follows from Lemma~13 in~\cite{Avd_solv} whose
proof is based on the classification of all possibilities for a pair
of active roots (see Theorem~\ref{thm_classification} below).
Further we give a proof of Proposition~\ref{prop_not_contained} that
does not use this classification.

\begin{proof}[Proof of
Proposition~\textup{\ref{prop_not_contained}}] We first introduce
the notation $I = \Supp \alpha \cup \Supp \beta$. In view of
Propositions~\ref{prop_bijection} and~\ref{prop_equal_pi}, for every
root $\gamma \in I$ we have a well-defined index $i(\gamma) \in \{1,
\ldots, m\}$ such that there exists a root $\widetilde \gamma \in
\Psi_{i(\gamma)}$ with $\pi(\widetilde \gamma) = \gamma$. We may
assume that $\widetilde \gamma \in F(\alpha) \cup F(\beta)$. If the
map $i \colon I \to \{1, \ldots, m\}$ is injective, then the set
$\tau(F(\alpha) \cup F(\beta))$ contains at least $|I|$ different
weights, which are linearly independent by
Theorem~\ref{thm_solv_spher}. This implies that all weights
$\tau(\gamma)$, where $\gamma \in I$, are linearly independent.
Hence we have $\tau(\alpha) \ne \tau(\beta)$, which is false.
Therefore there are two different roots $\gamma_1, \gamma_2 \in I$
with $i(\gamma_1) = i(\gamma_2)$. In view of
Propositions~\ref{prop_bijection} and~\ref{prop_F&Supp}(a,b,c) the
roots $\gamma_1, \gamma_2$ cannot lie simultaneously in the set
$\Supp \alpha$; by the same reason they cannot lie simultaneously in
the set $\Supp \beta$. Hence one of these two roots lies in $\Supp
\alpha \backslash \Supp \beta$ and the other one lies in $\Supp
\beta \backslash \Supp \alpha$. Without loss of generality we may
assume that $\gamma_1 \in \Supp \alpha \backslash \Supp \beta$ and
$\gamma_2 \in \Supp \beta \backslash \Supp \alpha$. Then the unique
root $\widetilde \alpha \in F(\alpha)$ with $\pi(\widetilde \alpha)
= \gamma_1$ and the unique root $\widetilde \beta \in F(\beta)$ with
$\pi(\widetilde \beta) = \gamma_2$ are the desired roots.
\end{proof}

\begin{theorem}[{\rm \cite[Theorem~3]{Avd_solv}}]\label{thm_active_roots}
For every root $\alpha \in \Psi$, the pair $(\alpha,\pi(\alpha))$ is
contained in Table~\textup{\ref{table_active_roots}}.
\end{theorem}

\begin{table}[h]

\caption{}\label{table_active_roots}

\begin{center}

\begin{tabular}{|c|c|c|c|}

\hline

No. & Type of $\Delta(\alpha)$ & $\alpha$ & $\pi(\alpha)$ \\

\hline

1 & any of rank $n$ & $\alpha_1 + \alpha_2 + \ldots + \alpha_n$ &
$\alpha_1, \alpha_2, \ldots, \alpha_n$\\

\hline

2 & $\mathsf B_n$ & $\alpha_1 + \alpha_2 + \ldots + \alpha_{n-1} +
2\alpha_n$ &
$\alpha_1, \alpha_2, \ldots, \alpha_{n-1}$\\

\hline

3 & $\mathsf C_n$ & $2\alpha_1 + 2\alpha_2 + \ldots + 2\alpha_{n-1}
+ \alpha_n$
& $\alpha_n$\\

\hline

4 & $\mathsf F_4$ & $2\alpha_1 + 2\alpha_2 + \alpha_3 + \alpha_4$ &
$\alpha_3, \alpha_4$\\

\hline

5 & $\mathsf G_2$ & $2\alpha_1 + \alpha_2$ & $\alpha_2$\\

\hline

6 & $\mathsf G_2$ & $3\alpha_1 + \alpha_2$ & $\alpha_2$\\

\hline
\end{tabular}

\end{center}

\end{table}

Let us explain the notation in Table~\ref{table_active_roots}. In
the column `$\alpha$' the expression of $\alpha$ as the sum of
simple roots in $\Supp \alpha$ is given. At that, the $i$th simple
root in the diagram $\Sigma(\Supp\alpha)$ is denoted by~$\alpha_i$.
In the column `$\pi(\alpha)$' we list all possibilities for
$\pi(\alpha)$ for a given active root~$\alpha$.

The following result is obtained by an easy case-by-case
consideration of all the possibilities in
Table~\ref{table_active_roots}.

\begin{corollary}[{\rm \cite[Corollary~7]{Avd_solv}}] \label{crl_terminal}
If $\alpha \in \Psi$, $|\Supp \alpha| \ge 2$, and ${\delta \in \Supp
\alpha \cap F(\alpha)}$, then $\delta$ is terminal with respect to
$\Supp \alpha$.
\end{corollary}

Further we give a list of several conditions on a pair $\alpha,
\beta$ of active roots. These conditions will be used below when we
formulate Theorem~\ref{thm_classification}.

\begin{wrapfigure}{R}{150\unitlength}

\begin{picture}(150,115)
\put(70,75){\circle{4}} \put(65,83){$\gamma_0$}

\put(70,73){\line(0,-1){18}} \put(70,53){\circle{4}}
\put(76,51){$\gamma_1$} \put(70,51){\line(0,-1){11}}
\put(70,36){\circle*{0.5}} \put(70,32){\circle*{0.5}}
\put(70,28){\circle*{0.5}} \put(70,25){\line(0,-1){11}}
\put(70,12){\circle{4}} \put(76,10){$\gamma_r$}

\put(72,76){\line(2,1){16}} \put(68,76){\line(-2,1){16}}

\put(90,85){\circle{4}} \put(88,74){$\beta_1$}
\put(92,86){\line(2,1){10}} \put(106,93){\circle*{0.5}}
\put(110,95){\circle*{0.5}} \put(114,97){\circle*{0.5}}
\put(118,99){\line(2,1){10}} \put(130,105){\circle{4}}
\put(128,94){$\beta_q$}

\put(50,85){\circle{4}} \put(45,74){$\alpha_1$}
\put(48,86){\line(-2,1){10}} \put(34,93){\circle*{0.5}}
\put(30,95){\circle*{0.5}} \put(26,97){\circle*{0.5}}
\put(22,99){\line(-2,1){10}} \put(10,105){\circle{4}}
\put(05,94){$\alpha_p$}

\end{picture}

\caption{}\label{diagram_difficult}

\end{wrapfigure}

$(\mathrm D0)$ $\Supp \alpha \cap \Supp \beta = \varnothing$;

$(\mathrm D1)$ $\Supp \alpha \cap \Supp \beta = \{\delta\}$, where
$\pi(\alpha) \ne \delta$, $\pi(\beta) \ne \delta$, and $\delta$ is
terminal with respect to both $\Supp \alpha$ and $\Supp \beta$;

$(\mathrm E1)$ $\Supp \alpha \cap \Supp \beta = \{\delta\}$, where
$\delta = \pi(\alpha) = \pi(\beta)$, $\alpha - \delta \in \Delta_+$,
$\beta - \delta \in \Delta_+$, and $\delta$ is terminal with respect
to both $\Supp \alpha$ and $\Supp \beta$;

$(\mathrm D2)$ the diagram $\Sigma(\Supp \alpha \cup \Supp \beta)$
has the form shown on Figure~\ref{diagram_difficult} (for some $p,
q, r\ge 1$), $\alpha = \alpha_1 + \ldots + \alpha_p + \gamma_0 +
\gamma_1 + \ldots + \gamma_r$, $\beta = \beta_1 + \ldots + \beta_q +
\gamma_0 + \gamma_1 + \ldots + \gamma_r$, $\pi(\alpha) \notin \Supp
\alpha \cap \Supp \beta$, and $\pi(\beta) \notin \Supp \alpha \cap
\Supp \beta$;

$(\mathrm E2)$ the diagram $\Sigma(\Supp \alpha \cup \Supp \beta)$
has the form shown on Figure~\ref{diagram_difficult} (for some
$p,q,r\ge 1$), $\alpha = \alpha_1 + \ldots + \alpha_p + \gamma_0 +
\gamma_1 + \ldots + \gamma_r$, $\beta = \beta_1 + \ldots + \beta_q +
\gamma_0 + \gamma_1 + \ldots + \gamma_r$, and $\pi(\alpha) =
\pi(\beta) \in \Supp \alpha \cap \Supp \beta$.

Now let us introduce the set $\Pi_0 = \bigcup\limits_{\delta \in
\mathrm M} \Supp \delta \subset \Pi$.

In addition to the set $\Upsilon_0(H) = (\mathrm M, \pi, \sim)$, to
the subgroup $H$ we also assign the set $\Upsilon(H) = {(S, \mathrm
M, \pi, \sim)}$, where $\pi$ is regarded as a map from $\mathrm M$
to~$\Pi$.

The following theorem gives a classification of all connected
solvable spherical subgroups in~$G$ standardly embedded in~$B$.

\begin{theorem} \label{thm_classification}
\textup{(a)~\cite[Theorem~4]{Avd_solv}} Up to conjugation by an
element in~$T$, $H$ is uniquely determined by the set~$\Upsilon(H) =
(S, \mathrm M, \pi, \sim)$, and this set satisfies the following
conditions:

$(\mathrm A)$ $\pi(\alpha) \in \Supp \alpha$ for every root $\alpha
\in \mathrm M$, and the pair $(\alpha, \pi(\alpha))$ is contained in
Table~\textup{\ref{table_active_roots}};

$(\mathrm D)$ if $\alpha, \beta \in \mathrm M$ and $\alpha \nsim
\beta$, then one of possibilities $(\mathrm D0)$, $(\mathrm D1)$,
$(\mathrm D2)$ is realized;

$(\mathrm E)$ if $\alpha, \beta \in \mathrm M$ and $\alpha \sim
\beta$, then one of possibilities $(\mathrm D0)$, $(\mathrm D1)$,
$(\mathrm E1)$, $(\mathrm D2)$, $(\mathrm E2)$ is realized;

$(\mathrm C)$ if $\alpha \in \mathrm M$, then $\Supp \alpha
\not\subset \bigcup\limits_{\delta \in \mathrm M \backslash
\{\alpha\}}\Supp \delta$;

$(\mathrm T)$ $\langle \Ker \tau \rangle \cap \langle \Pi_0 \rangle
= \langle \alpha - \beta \mid \alpha, \beta \in \mathrm M, \alpha
\sim \beta \rangle$.

\textup{(b)~\cite[Theorem~5]{Avd_solv}} Suppose that a subtorus $S
\subset T$, a subset $\mathrm M \subset \Delta_+$, a map $\pi \colon
\mathrm M \to \Pi$, and an equivalence relation $\sim$ on $\mathrm
M$ satisfy conditions $(\mathrm A)$, $(\mathrm D)$, $(\mathrm E)$,
$(\mathrm C)$, and $(\mathrm T)$. Then there exists a connected
solvable spherical subgroup $H \subset G$ standardly embedded in~$B$
such that $\Upsilon(H) = (S, \mathrm M, \pi, \sim)$.
\end{theorem}

\begin{remark}
Generally speaking, it may happen that $\Upsilon(H) \ne
\Upsilon(H')$, but the subgroups $H, H'$ are conjugate in~$G$.
Concerning the problem of when two connected solvable spherical
subgroups standardly embedded in~$B$ are conjugate in~$G$, see
Proposition~\ref{prop_two_conj} and Theorem~\ref{thm_elem_tr} below
(a more detailed discussion of this problem see also
in~\cite[\S\,5]{Avd_solv}).
\end{remark}

\begin{proposition}[{\rm\cite[Proposition~13]{Avd_solv}}] \label{prop_two_conj}
Suppose that $H$, $H'$ are such that $H' = g H g^{-1}$ for some $g
\in G$. Then ${g \in N' \cdot N_G(T) \cdot N}$.
\end{proposition}

\begin{definition}
Suppose that $\delta \in \Psi^{\reg} \cap \Pi$. An
\textit{elementary transformation with center~$\delta$} (or simply
an \textit{elementary transformation}) is a transformation $H
\mapsto \rho_\delta H \rho_\delta^{-1}$, where $\rho_\delta \in
\theta^{-1}(r_\delta)$.
\end{definition}

As can be easily seen, in this definition the subgroup $\rho_\delta
H \rho_\delta^{-1}$ is also standardly embedded in~$B$.

The argument used in the proof of Theorem~6 in~\cite{Avd_solv}
actually proves the following theorem.

\begin{theorem} \label{thm_elem_tr}
Suppose that $H$, $H'$ are such that $H' = \sigma H \sigma^{-1}$ for
some ${\sigma \in N_G(T)}$. Then there are elements $\rho_1, \ldots,
\rho_k \in N_G(T)$ with the following properties:

\textup{(1)} $\theta(\rho_1), \ldots, \theta(\rho_k)$ are simple
reflections;

\textup{(2)} $\sigma = \rho_k \ldots \rho_1$;

\textup{(3)} the chain $H \mapsto \rho_1 H\rho_1^{-1} \mapsto
\rho_2\rho_1 H \rho_1^{-1}\rho_2^{-1} \mapsto \ldots \mapsto \sigma
H \sigma^{-1} = H'$ is a chain of elementary transformations.
\end{theorem}

\begin{proposition}[{\rm\cite[Lemma~30(a,b)]{Avd_solv}}]
\label{prop_changes_Psi&pi} Suppose that $H \mapsto H'$ is an
elementary transformation with center~$\delta$. Then:

\textup{(a)} $\Psi' = r_\delta(\Psi \backslash \{\delta\}) \cup
\{\delta\}$;

\textup{(b)} $\pi'(r_\delta(\alpha)) = \pi(\alpha)$ for $\alpha \in
\Psi \backslash \{\delta\}$;
\end{proposition}

\begin{corollary} \label{crl_etar}
Suppose that $H \mapsto H'$ is an elementary transformation with
center~$\delta$ and $\alpha \in \Psi$ is an arbitrary root. Then for
the unique positive root $\widetilde \alpha$ in the set
$\{r_\delta(\alpha), -r_\delta(\alpha)\}$ we have $\widetilde \alpha
\in \Psi'$ and $\pi'(\widetilde \alpha) = \pi(\alpha)$.
\end{corollary}

Theorem~\ref{thm_elem_tr} and Corollary~\ref{crl_etar} imply the
following corollary.

\begin{corollary} \label{crl_assoc_coinc}
Suppose that $\sigma \in N_G(H) \cap N_G(T)$ and $\alpha \in \Psi$
is an arbitrary root. Put $s = \theta(\sigma) \in W$. Then for the
unique positive root~$\beta$ in the set $\{s(\alpha), -s(\alpha)\}$
we have $\beta \in \Psi$ and $\pi(\beta) = \pi(\alpha)$.
\end{corollary}

\section{Proof of the main results} \label{sect_proof}

In this subsection we prove Theorem~\ref{thm_main_result} and
Proposition~\ref{prop_root_in_P}. We preserve the notation
introduced in
\S\,\ref{sect_preliminaries}--\ref{sect_further_results}.

\begin{proof}[Proof of Theorem~\textup{\ref{thm_main_result}}]
In view of Proposition~\ref{prop_two_conj} we have $N_G(H) \subset N
\cdot N_G(T) \cdot N$. Since $N \subset H \subset N_G(H)$, it
remains to find the group $N_G(H) \cap N_G(T)$.

Suppose that $\sigma \in N_G(T)$ is such that $\sigma H \sigma^{-1}
= H$. Put $s = \theta(\sigma) \in W$. Theorem~\ref{thm_elem_tr} and
Proposition~\ref{prop_changes_Psi&pi} imply that $s$ is contained in
the subgroup~$W_0$ of~$W$ generated by simple reflections
corresponding to the roots in~$\Pi_0$.

If $s = e$, then $\sigma \in N_G(H) \cap T$. It was shown in
\S\,\ref{sect_result} that $N_G(H) \cap T = A$.

In what follows we assume $s \ne e$.

\begin{lemma} \label{lemma_positive}
If $\alpha \in \Delta_+ \backslash \Psi^{\reg}$, then $s(\alpha) \in
\Delta_+$.
\end{lemma}

\begin{proof}
From the definition of the set $\Psi^{\reg}$ it follows that the
projection of the subspace~$\mathfrak h$ to $\mathfrak g_\alpha$ is
not zero. Since $(\Ad \sigma)\mathfrak g_\gamma = \mathfrak
g_{s(\gamma)}$ for every $\gamma \in \Delta$, the condition $(\Ad
\sigma)\mathfrak h = \mathfrak h$ implies that the projection of
$\mathfrak h$ to $\mathfrak g_{s(\alpha)}$ is not zero as well.
Hence $s(\alpha) \in \Delta_+$.
\end{proof}

\begin{proposition} \label{prop_two_possib}
For every $\alpha \in \Psi$ we have $s(\alpha) \in \{\alpha, -
\alpha\}$.
\end{proposition}

\begin{proof}
Let $\widetilde \alpha$ be the unique positive root in the set
$\{s(\alpha), -s(\alpha)\}$. Our aim is to prove that $\widetilde
\alpha = \alpha$. In view of Corollary~\ref{crl_assoc_coinc} we have
${\widetilde \alpha \in \Psi}$ and $\pi(\widetilde \alpha) =
\pi(\alpha)$. If the set $\pi^{-1}(\pi(\alpha))$ consists of one
element, then $\widetilde \alpha = \alpha$ and the assertion is
proved. Further we assume that the set $\pi^{-1}(\pi(\alpha))$
contains more than one element. Let $\beta \in \pi^{-1}(\pi(\alpha))
\backslash \{\alpha\}$ be an arbitrary root. Then by
Proposition~\ref{prop_equal_pi} we have $\tau(\alpha) =
\tau(\beta)$. Further, in view of
Proposition~\ref{prop_not_contained} there are roots $\alpha' \in
F(\alpha)$ and $\beta' \in F(\beta)$ such that $\tau(\alpha') =
\tau(\beta')$ and ${\pi(\alpha') \notin \Supp \beta}$. All the four
roots $\alpha$, $\beta$, $\alpha'$, $\beta'$ are active and pairwise
different, therefore none of them is regular. In particular, each of
the roots $\alpha$, $\alpha'$ is contained in the set~$\Delta_+
\backslash \Psi^{\reg}$. Besides, for the root $\alpha'' = \alpha -
\alpha' \in \Delta_+$ we have $\alpha'' \in \Delta_+ \backslash
\Psi$, whence by Lemma~\ref{lemma_positive} all the three roots
$s(\alpha)$, $s(\alpha')$, and $s(\alpha'')$ are positive. In
particular, $\widetilde \alpha = s(\alpha)$. Applying
Corollary~\ref{crl_assoc_coinc}, we obtain $s(\alpha), s(\alpha')
\in \Psi$ and $\pi(s(\alpha')) = \pi(\alpha')$. Now the equality
$s(\alpha) = s(\alpha') + s(\alpha'')$ implies that $\pi(\alpha')
\in \Supp s(\alpha)$. The latter means that $\widetilde \alpha \ne
\beta$, which completes the proof.
\end{proof}

Proposition~\ref{prop_two_possib} and Lemma~\ref{lemma_positive}
imply the following corollary.

\begin{corollary} \label{crl_not_reg}
If a root $\alpha \in \Psi$ is not regular, then $s(\alpha) =
\alpha$.
\end{corollary}

\begin{proposition} \label{prop_height_ge2}
Suppose that $\alpha \in \Psi$ and $\hgt \alpha \ge 2$. Then
$s(\alpha) = \alpha$.
\end{proposition}

\begin{proof}
The condition $\hgt \alpha \ge 2$ implies that $|\Supp \alpha| \ge
2$, whence ${|F(\alpha)| \ge 2}$ by
Proposition~\ref{prop_F&Supp}(a). Let $\alpha' \in F(\alpha)
\backslash \{\alpha\}$ be an arbitrary root. Since $\alpha - \alpha'
\in {\Delta_+ \backslash \Psi \subset \Delta_+ \backslash
\Psi^{\reg}}$, by Lemma~\ref{lemma_positive} we obtain that
$s(\alpha) - s(\alpha') \in \Delta_+$. Further, in view of
Proposition~\ref{prop_two_possib} we have ${s(\alpha) \in \{\alpha,
-\alpha\}}$ and ${s(\alpha') \in \{\alpha', -\alpha'\}}$. If
$s(\alpha) = -\alpha$, then the root $s(\alpha) - s(\alpha')$ lies
in the set $\{{-\alpha + \alpha'}, {-\alpha - \alpha'}\}$, which
evidently does not contain positive roots. Hence $s(\alpha) =
\alpha$.
\end{proof}

In view of Corollary~\ref{crl_not_reg},
Proposition~\ref{prop_height_ge2} implies the following corollary.

\begin{corollary} \label{crl_simple_reg}
If a root $\alpha \in \Psi$ has the property $s(\alpha) = - \alpha$,
then $\alpha \in {\Psi^{\reg} \cap \Pi}$.
\end{corollary}

\begin{proposition} \label{prop_orthogonal}
Suppose that a root $\alpha \in \Psi^{\reg} \cap \Pi$ has the
property $s(\alpha) = -\alpha$. Then $(\alpha, \beta) = 0$ for every
root $\beta \in \Psi \backslash \{\alpha\}$.
\end{proposition}

\begin{proof}
Since the group $W$ acts on the space $Q$ by orthogonal
transformations, we have $(\alpha, \beta) = (s(\alpha), s(\beta))$.
Assume that ${(\alpha, \beta) \ne 0}$. Then the hypothesis and
Proposition~\ref{prop_two_possib} imply that $s(\beta) = -\beta$. By
Corollary~\ref{crl_simple_reg} we obtain that ${\beta \in \Pi}$,
whence $(\alpha, \beta) < 0$. Therefore $\alpha + \beta \in
\Delta_+$ and $s(\alpha + \beta) = s(\alpha) + s(\beta) = {-(\alpha
+ \beta)} \notin \Delta_+$. On the other hand, it is easy to see
that $\alpha + \beta \notin \Psi$, whence by
Lemma~\ref{lemma_positive} we obtain that ${s(\alpha + \beta) \in
\Delta_+}$, a contradiction.
\end{proof}

Corollary~\ref{crl_simple_reg} and Proposition~\ref{prop_orthogonal}
imply the following corollary.

\begin{corollary} \label{crl_root_in_P}
Every root $\alpha \in \Psi$ with $s(\alpha) = - \alpha$ is
contained in~$\mathrm P$.
\end{corollary}

Let $\delta_1, \ldots, \delta_k$ be all the roots in the set
$\{\alpha \in \Psi \mid s(\alpha) = -\alpha\}$. By
Corollary~\ref{crl_root_in_P} all of them are contained in~$\mathrm
P$. Next, it is easy to see that every active root is invariant
under the action of the element $r_{\delta_1} \ldots r_{\delta_k} s
\in W_0$. It follows from Proposition~\ref{prop_F&Supp}(c) that
$\Pi_0 \subset \langle \Psi \rangle$, whence $r_{\delta_1} \ldots
r_{\delta_k} s = e$. Therefore $s = r_{\delta_k} \ldots
r_{\delta_1}$ and $\sigma = t \rho_{\delta_k} \ldots
\rho_{\delta_1}$ for some $t \in A$.

Each of the roots $\delta_1, \ldots, \delta_k$ is contained in the
space $\langle \Pi_0 \rangle$ and is orthogonal to the subspace
$\langle L \rangle$. Further, from condition~$(\mathrm T)$ (see
Theorem~\ref{thm_classification}) we have $\langle \Ker \tau \rangle
\cap \langle \Pi_0 \rangle = \langle L \rangle$. In view of what we
have said above the condition $s\langle \Ker \tau \rangle = \langle
\Ker \tau \rangle$ implies that the space $\langle \Ker \tau
\rangle$ is orthogonal to each of the roots $\delta_i$, $i = 1,
\ldots, k$. The latter means that $\delta_1, \ldots, \delta_k \in
\mathrm P_S$.

Thus we have obtained that $\sigma = t \rho_{\delta_k} \ldots
\rho_{\delta_1}$ for some roots $\delta_1, \ldots, \delta_k \in
\mathrm P_S$ and some element $t \in A$. On the other hand, as can
be easily seen, for every root $\delta \in \mathrm P_S$ we have
$\rho_\delta \in N_G(H)$. This completes the proof of the theorem.
\end{proof}

\begin{proof}[Proof of Proposition~\textup{\ref{prop_root_in_P}}]
Suppose that $\alpha \in \Pi$. First assume that $\alpha \in \mathrm
P$. Then there are two possibilities.

$1^\circ$. $\alpha \in \mathrm M$. Then it is easy to check that all
the conditions of Case~1 are fulfilled.

$2^\circ$. $\alpha \notin \mathrm M$. Then $\alpha \in \Supp \beta$
for some root $\beta \in \mathrm M$. In view of
Proposition~\ref{prop_F&Supp}(d) we have $\alpha \in F(\beta)$, and
Corollary~\ref{crl_terminal} implies that $\alpha$ is terminal with
respect to $\Supp \beta$. Having performed a case-by-case
consideration of all the possibilities in
Table~\ref{table_active_roots}, we find that the condition $(\alpha,
\beta) = 0$ holds if and only if the diagram $\Sigma(\Supp \beta)$
is of type $\mathsf B_n$ ($n \ge 2$) and conditions (1)--(4) of
Case~2 are fulfilled. Since a node of a Dynkin diagram cannot be
incident to two double edges, in view of conditions $(\mathrm D)$
and $(\mathrm E)$ (see Theorem~\ref{thm_classification}) for every
root $\beta' \in \mathrm M \backslash \{\beta\}$ we have $\alpha
\notin \Supp \beta'$, hence condition~(5) of Case~2.

Now let us prove the converse implication of the proposition.

If Case~1 takes place, then in view of Corollary~\ref{crl_either_or}
we obtain $\alpha \in \mathrm P$.

Suppose that Case~2 takes place. From conditions~(1)--(3) it follows
that $\alpha \in \Psi$. Conditions (1), (4), and~(5) imply that
$(\alpha, \delta) = 0$ for every root $\delta \in \Psi \backslash
\{\alpha\}$. It remains to prove that $\alpha \in \Psi^{\reg}$. If
this is not the case then there is a root $\alpha' \in \Psi
\backslash \{\alpha\}$ such that $\tau(\alpha) = \tau(\alpha')$. In
view of Proposition~\ref{prop_Psi's} and
Corollary~\ref{crl_either_or} this implies that $\beta' = \alpha' +
{(\beta - \alpha) \in \mathrm M}$. Then the diagram
$\Sigma(\{\alpha\} \cup \Supp \beta')$ is connected, and we have a
contradiction with condition~(5). Thus, $\alpha \in \mathrm P$.
\end{proof}


\end{document}